\documentclass[12pt]{amsart}

\usepackage{amsmath}
\usepackage{amsfonts}
\usepackage{amssymb,enumerate}
\usepackage{amsthm}

\usepackage{hyperref}
\usepackage{cleveref}
\usepackage{xcolor}
\usepackage{enumitem}
\usepackage{float}
\usepackage{verbatim}
\usepackage{tikz}
\usetikzlibrary{arrows}
\usepackage{mathrsfs}
\usetikzlibrary{arrows}
\usepackage{listings}
\definecolor{DimGray}{rgb}{0.41, 0.41, 0.41}

\theoremstyle{definition}
\newtheorem{defin}{Definition}
\newtheorem{question}{Question}

\theoremstyle{plain}

\newtheorem{theorem}{Theorem}
\newtheorem*{theorem*}{Theorem}

\theoremstyle{remark}

\title{The 4/3 problem for germs of isolated plane curve singularities}

\author{Patricio Almir\'on}

\address{Instituto de Matemática Interdisciplinar (IMI) and Departamento de \'{A}lgebra, Geometr\'{i}a y Topolog\'{i}a\\
Facultad de Ciencias Matem\'{a}ticas,\\
Universidad Complutense de Madrid\\
28040, Madrid, Spain.}


\thanks{The author is supported by Spanish Ministerio de Ciencia, Innovaci\'{o}n y Universidades MTM2016-76868-C2-1-P}

\email{palmiron@ucm.es}

\begin{document}

\begin{abstract}
In this survey we overview the different approaches and solutions of a question posed by Dimca and Greuel about the quotient of the Milnor and Tjurina numbers.
\end{abstract}

\maketitle

\section{Introduction}
Analytic and topological invariants of germs of isolated plane curve singularities are central objects in Singularity Theory, see \cite{greuel-book} and the references therein for an overview. One of the main objects of study is to find relations between them and to find topological constrains for analytical invariants. As one can see in \cite{dim} and \cite{greuel-book}, two mainstream invariants are the Milnor number \(\mu,\) and the Tjurina number \(\tau.\) In fact, the Milnor number is a topological invariant and the Tjurina number an analytic invariant. If \(C:=\{f(x,y)=0\}\) is a germ of isolated plane curve singularity, the easiest way to define these numbers is:
\[ \mu := \dim_{\mathbb{C}} \frac{ \mathbb{C} \{x, y\} }{ (\partial f/ \partial x, \partial f/ \partial y) }, \quad \tau := \dim_{\mathbb{C}} \frac{ \mathbb{C} \{x, y\} }{ (f, \partial f/ \partial x, \partial f/ \partial y) }. \]

In 2017 on \cite{dim}, Dimca and Greuel posed the following question:
\begin{question} \label{conjecture}
Is it true that $\mu/\tau < 4/3$ for any isolated plane curve singularity?
\end{question}
This guessed bound is inferred by Dimca and Greuel from some families of plane curves that asymptotically achieve this bound. From this point view, Question \ref{conjecture} can be divided into two questions: is it true? If it is true, can  the \(4/3\) bound be inferred from the geometry of the plane curve singularity? 

In this survey, we will try to show the different approaches to this question and the main problems attached to them. 
 \section{Deformation theory}
In this section, we are going to show that Milnor and Tjurina numbers are closely related to the theory of deformations. We refer to \cite{greuel-book} for general deformation theory.

\begin{defin}
Let \((C,0)\) be a germ of  isolated plane curve singularity. A deformation of \((C,0)\) is a germ of flat morphism \((\mathcal{Y},0)\rightarrow(S,0)\) whose special fibre is isomorphic to \((C,0)\). We call \((S,0)\) the base space of the deformation. The deformation is called versal if any other deformation results from it by base change. It is called miniversal if it is versal and \(S\) has minimal possible dimension. 
\end{defin}
In \cite{greuel-book}, it is shown that an explicit way to construct versal and miniversal deformations of a plane curve is by using the Milnor algebra \(M_f,\) and the Tjurina algebra \(T_f.\)
\[ M_f := \frac{ \mathbb{C} \{x, y\} }{ (\partial f/ \partial x, \partial f/ \partial y) }, \quad T_f := \frac{ \mathbb{C} \{x, y\} }{ (f, \partial f/ \partial x, \partial f/ \partial y) }. \]
\begin{theorem}[(Tjurina) Corollary 1.17 \cite{greuel-book}]
Let \((C,0)\) be a germ of  isolated plane curve singularity defined by \(f\in\mathcal{O}_{\mathbb{C}^2,0}\) and \(g_1,\dots,g_k\in\mathcal{O}_{\mathbb{C}^2,0}\) be a \(\mathbb{C}\)--basis of \(T_f\) (resp. of \(M_f\))
If we set,
\[F(x,\mathbf{t}):=f(x)+\sum_{j=1}^{k}t_jg_j(x),\quad (\mathcal{X},0):=V(F)\subset(\mathbb{C}^2\times\mathbb{C}^k,0),\]
then \((C,0)\hookrightarrow(\mathcal{X},0)\xrightarrow{\varphi}(\mathbb{C}^k,0),\) with \(\varphi\) the projection from the second component, is a miniversal (resp. versal) deformation of \((C,0).\)
\end{theorem}

Inside the base space of a miniversal deformation of a germ of plane curve singularity there is an interesting closed analytic subspace \(\Delta^{\mu}\) called the \(\mu\)--constant stratum. This stratum can be defined as follows: take the miniversal deformation \(\varphi:(\mathcal{Y},0)\rightarrow(S,0)\) of a plane curve singularity \(C.\) Denote by \(\mu\) the Milnor number of \(C\) and by \(\mathcal{Y}_s:=\varphi^{-1}(s)\) a fiber of the deformation, then 
\[\Delta^{\mu}:=\{s\in S\;|\;\mu(\mathcal{Y}_s)=\mu\}.\] 
Then it can be proven that this stratum is smooth (see Theorem 2.61 in \cite{greuel-book}) and its codimension can be computed from the embedded resolution of the plane curve by the following formula given by Wall in Section 8 of \cite{wall2}.
\begin{theorem}[Theorem 8.1 in \cite{wall2}, (2.8.36) pg. 373 in \cite{greuel-book} ] \label{codim}
If \((C,0)\) is a germ of plane curve singularity, \(e_p\) is the sequence of multiplicities of the strict transform of the embedded resolution of \(C\) and \(c\) is the number of free points in the resolution then 
\[\operatorname{codim}(\Delta^{\mu})=\sum_p\frac{e_p(e_p+1)}{2}-c-1.\]
\end{theorem}

\section{Solutions to Dimca and Greuel question}
Dimca and Greuel's question has been completely solved by the author in \cite{quotpat}. However, before this general solution, there has been several solutions from different points of view for some families of plane curve singularities. In this section, we will try to overview in a chronological order the different results until reaching the general solution of Dimca and Greuel's question.

The first result about Question \ref{conjecture} was given in 2018 by Blanco and the author in \cite{alblanc} for semi-quasi-homogeneous singularities. We recall that \( f \) is a semi-quasi-homogeneous singularity with weights \( w = (n, m) \) such that \( \gcd(n, m) \geq 1 \) and \( n, m \geq 2 \) if  \( f = f_0 + g \) is a deformation of the initial term \( f_0 = y^n - x^m \) such that \( \deg_w(f_0) < \deg_w(g) \). For such singularities, Blanco and the author in \cite{alblanc} give a positive answer to this question. This answer is due to a formula for the minimal Tjurina number of the family of semi-quasi-homogeneous given by Briançon, Granger and Maisonobe in \cite{granger}. The idea here is to use the upper semicontinuity of the Tjurina number (Theorem 2.6 in \cite{greuel-book}) to reduce the proof to show the inequality for \(\mu/\tau_{\min}.\) In fact, until the appearance of the general solution this was the only non--irreducible family of plane curve singularities for which Question \ref{conjecture} was solved.

In 2019, a series of three preprints \cite{taumin}, \cite{genzmertau}, \cite{wang} appeared in a short time. They give a positive answer for the case of irreducible germs of plane curve singularities. The three approaches are based on the explicit computation of the dimension of the generic component of the moduli space of irreducible plane curve singularities given by Genzmer in \cite{genzmer16} in terms of the sequence of multiplicities of the strict transform of a resolution of the irreducible plane curve. It was shown by Zariski in \cite{zariski-moduli} and Teissier in the appendix to the book of Zariski \cite{teissier-appendix} that to compute the dimension of the generic component of the moduli space of irreducible plane curves is closely related to compute the minimal Tjurina number in the equisingularity class of a branch. The relation of the dimension of the generic component of the moduli and the minimal Tjurina number of irreducible plane curves is due to the properties of Teissier monomial curve (see \cite{teissier-appendix}).

In April 2019, Alberich-Carrami\~{n}ana, Blanco, Melle-Hern\'{a}ndez and the author in \cite{taumin} gave a positive answer to Question \ref{conjecture} through a formula for the minimal Tjurina number in an equisingularity class of irreducible plane curve singularities in terms of the sequence of multiplicities that can be obtained from Genzmer's formula in \cite{genzmer16} together with Wall's formula (Theorem \ref{codim}). A few days after, Genzmer and Hernandes in \cite{genzmertau} provided an alternative proof of Dimca and Greuel's inequality. Even if the techniques used are quite different, both results are based on the explicit computations given by Genzmer in \cite{genzmer16}. Finally, at the end of April, Wang in \cite{wang} gave another alternative proof for the irreducible case based also in Genzmer's result about the dimension of the generic component of the moduli space in \cite{genzmer16}. However, Wang's approach is very interesting since he proves that \(3\mu-4\tau\) is a monotonic increasing invariant under blow--ups for irreducible plane curve singularities which provides a nice perspective in the possible applications of Dimca and Greuel's question.

After the previous discussion one realized that all the answers to Dimca and Greuel's question are based on the explicit computation of the minimal Tjurina number of certain families of singularities. However, they do not provide an answer to the second question that we formulated in the introduction: can the \(4/3\) bound be inferred from the geometry of the plane curve singularity?
One may think that Wang's approach gives the answer to this question for the irreducible case. However, one cannot prove that \(3\mu-4\tau\) is an increasing monotonic invariant under blow-ups if one does not have Genzmer's formula for the dimension of the generic component of the moduli space. Moreover, implicitly there is no reason to consider \(a\mu-b\tau\) with \((a,b)\neq(3,4).\) In this way, the reason about \(4/3\) remained open after Wang's result even for irreducible plane curves.

Finally, based on the idea to provide a full answer to Dimca and Greuel's question, the author in \cite{quotpat} changed the point of view. This new approach is based on the theory of deformation for surface singularities. More concretely, in the geometry of normal two--dimension double point singularities. Normal two--dimension double point singularities have equation \(\{z^2+f(x,y)=0\}\subset\mathbb{C}^3\) with \(f(x,y)=0\) defining a germ of plane curve singularity. Moreover, they have the same Milnor and Tjurina numbers than the associated plane curve singularity. From this point of view, one can use the upper bound for the difference \(\mu-\tau\) given by Wahl in \cite{wahl} in terms of the geometric genus of the double point singularity. After that, the good properties of the geometric genus of such a surface singularities allow the author to provide a full answer to Dimca and Greuel's question.
\begin{theorem}[\cite{quotpat}]\label{teorema4/3}
For any germ of plane curve singularity 
\[\frac{\mu}{\tau}<\frac{4}{3}.\]
\end{theorem}
Moreover, in this setting the bound \(4/3\) is inferred from the geometric genus of the surface singularity. Also, the author in \cite{quotpat} provides a general framework that allow to continue with the problem of finding bounds for the quotient of Milnor and Tjurina numbers in higher dimension.


\end{document}